\documentclass[12pt]{article} 
\usepackage{amsmath}
\usepackage{amssymb} 
\usepackage{amscd}
\usepackage{amsthm}
\usepackage[utf8]{inputenc}
\usepackage{stmaryrd}
\usepackage{hyperref}
\newtheorem{theorem}{Theorem}

\newtheorem{proposition}[theorem]{Proposition}

\def\gr{\mathrm{gr}}

\def\dim{{\mbox{dim}}}

\def\ker{{\mbox{Ker}}}

\def\Hom {{\mbox{Hom}}}
\def\Ext {{\mbox{Ext}}}

\def\End{{\mbox{End}}}

\def\cala{{\cal A}}

\def\fraca{{\mathfrak A}}
 
 \def\fracg{{\mathfrak g}}

\def\bbbone{\mbox{\rm 1\hspace {-.6em} l}}

\def\Tor{{\mbox{Tor}}}

\numberwithin{equation}{section}

\begin{document}
\enlargethispage{3cm}

\thispagestyle{empty}
\begin{center}
{\bf FROM KOSZUL DUALITY TO POINCAR\'E DUALITY}
\end{center} 
\vspace{0.3cm}

\begin{center}
 Michel DUBOIS-VIOLETTE
\footnote{Laboratoire de Physique Th\'eorique, CNRS UMR 8627\\ 
Universit\'e Paris Sud 11,
B\^atiment 210\\ F-91 405 Orsay Cedex\\
Michel.Dubois-Violette$@$u-psud.fr}\\
\end{center}

 \vspace{0,5cm}
 
 \begin{center}
 {\sl Dédié à Raymond Stora}
 \end{center}

 \vspace{0,5cm}
\begin{abstract}
We discuss the notion of Poincaré duality for graded algebras and its connections with the Koszul duality for quadratic Koszul algebras. The relevance of the Poincaré duality is pointed out for the existence of twisted potentials associated to Koszul algebras as well as for the extraction of a good generalization of  Lie algebras among the quadratic-linear algebras.
\end{abstract}
\vfill

\noindent LPT-ORSAY 11-78
\baselineskip=0,7cm
\section*{Introduction}
Koszul complexes and Koszul resolutions play a major role at several places in theoretical physics. In many cases, this is connected with the BRST methods introduced in \cite{bec-rou-sto:1974}, \cite{bec-rou-sto:1975}, \cite{bec-rou-sto:1976a}, \cite{bec-rou-sto:1976b}, \cite{tyu:1975}.
For instance they enter in the classical BRST approach to constrained systems, in particular to gauge theory, and they are involved in the renormalization process of quantum field theory, see e.g.  \cite{mcm:1984}, \cite{mdv:1987c}, \cite{fis-hen-sta-tei:1989}, \cite{hen-tei:1992}, \cite{sto:2005}.\\

Our aim here is to describe elements of the formulation of the Koszul duality and of the Poincaré duality for quadratic algebras and for non homogeneous quadratic algebras and to draw some important consequence of the Poincaré duality property for Koszul algebras.\\

Throughout this paper $\mathbb K$ denotes a (commutative) field and all vector spaces, algebras, etc. are over $\mathbb K$. We use everywhere the Einstein summation convention over the repeated up-down indices.

\section{Graded algebras and Poincaré duality}
\subsection{Graded algebras}

The class of graded algebras involved here is the class of connected graded algebras $\cala$ of the form $\cala=T(E)/I$ where $E$ is a finite-dimensional vector space and where $I$ is a finitely generated graded ideal of the tensor algebra $T(E)$ such that $I=\oplus_{n\geq 2}I_n\subset \oplus_{n\geq 2}E^{\otimes^n}$. This class together with the homomorphisms of graded-algebras (of degree 0) as morphisms define the category $\mathbf{GrAlg}$.\\

For such an algebra $\cala=T(E)/I\in \mathbf{GrAlg}$ choosing a basis $(x^\lambda)_{\lambda\in\{1,\dots,d\}}$ of $E$ and a system of homogeneous independent generators $(f_\alpha)_{\alpha\in \{1,\dots,r\}}$ of $I$ with $(f_\alpha)\in E^{\otimes^{N_\alpha}}$ and $N_\alpha\geq 2$ for $\alpha\in \{1,\dots,r\}$, one can also write
\[
\cala=\mathbb K\langle x^1,\dots, x^d\rangle/[f_1,\dots, f_r]
\]
where $[f_1,\dots,f_r]$ is the ideal $I$ generated by the $f_\alpha$. Define $M_{\alpha\lambda}\in E^{\otimes^{N_\alpha-1}}$ by setting $f_\alpha=M_{\alpha\lambda}\otimes x^\lambda\in E^{\otimes^{N_\alpha}}$. Then the presentation of $\cala$ by generators and relations is equivalent to the exactness of the sequence of left $\cala$-modules
\begin{equation}
\cala^r \stackrel{M}{\rightarrow} \cala^d \stackrel{x}{\rightarrow} \cala \stackrel{\varepsilon}{\rightarrow} \mathbb K \rightarrow 0
\label{pres}
\end{equation}
where $M$ means right multiplication by the matrix $(M_{\alpha\lambda})\in M_{d,r}(\cala)$, $x$ means right multiplication by the column $(x^\lambda)$ and where $\varepsilon$ is the projection onto $\cala_0=\mathbb K$,  \cite{art-sch:1987}. In more intrinsic notations the exact sequence (\ref{pres}) reads
\begin{equation}
\cala\otimes R\rightarrow \cala\otimes E \stackrel{m}{\rightarrow} \cala \stackrel{\varepsilon}{\rightarrow} \mathbb K \rightarrow 0
\label{Ipres}
\end{equation}
where $R$ is the graded subspace of $T(E)$ spanned by the $f_\alpha$ $(\alpha\in\{1,\dots,r\})$, $m$ is the product in $\cala$ (remind that $E=\cala_1$) and where the first arrow is as in (\ref{pres}).\\

When $R$ is homogeneous of degree $N$ ($N\geq 2$), i.e. $R\subset E^{\otimes^N}$,  then $\cala$ is said to be a $N$-{\sl homogeneous algebra}: for $N=2$ one speaks of a quadratic algebra, for $N=3$ one speaks of a cubic algebra, etc.

\subsection{Global dimension}
 The exact sequence (\ref{Ipres}) of presentation of $\cala$ can be extended as a minimal projective resolution of the trivial left module $\mathbb K$, i.e. as an exact sequence of left modules
 \[
 \cdots \rightarrow M_n \rightarrow \cdots \rightarrow M_2 \rightarrow M_1 \rightarrow M_0 \rightarrow \mathbb K \rightarrow 0
 \]
 where the $M_n$ are projective i.e. in this graded case free left-modules \cite{car:1958}, which is minimal ; one has $M_0=\cala$, $M_1=\cala\otimes E$, $M_2=\cala\otimes R$ and more generally here $M_n =\cala\otimes E_n$ where the $E_n$ are finite-dimensional vector spaces. If such a minimal resolution has finite length $D<\infty$, i.e. reads 
 \begin{equation}
0\rightarrow \cala\otimes E_D \rightarrow \cdots \rightarrow \cala\otimes E \rightarrow \cala\rightarrow \mathbb K\rightarrow 0
\label{Mres}
\end{equation}
with $E_D\not=0$, then $D$ is an invariant called the {\sl left projective dimension} of $\mathbb K$ and it turns out that $D$ which coincide with the right projective dimension of $\mathbb K$ is also the sup of the lengths of the minimal projective resolutions of the left and of the right $\cala$-modules \cite{car:1958} which is called the {\sl global dimension} of $\cala$. Furthermore it was recently shown \cite{ber:2005} that this global dimension $D$ also coincides with the Hochschild dimension in homology as well as in cohomology. Thus for an algebra $\cala\in \mathbf{GrAlg}$, there is a unique notion of dimension from a homological point of view which is its global dimension $g\ell\dim(\cala)=D$ whenever it is finite.

\subsection{Poincaré duality versus Gorenstein property}

Let $\cala\in \mathbf{GrAlg}$ be of finite global dimension $D$. Then one has a minimal free resolution
\[
0\rightarrow M_D\rightarrow \cdots\rightarrow M_0\rightarrow \mathbb K \rightarrow 0
\]
 with $M_n=\cala\otimes E_n$, $\dim(E_n)<\infty$ and $E_2\simeq R$, $E_1\simeq E$ and $E_0\simeq \mathbb K$. By applying the functor $\Hom_\cala(\bullet, \cala)$ to the chain complex of free left $\cala$-module
 \begin{equation}
0\rightarrow M_D\rightarrow \cdots \rightarrow M_0\rightarrow 0
\label{M}
\end{equation}
one obtains the cochain complex
\begin{equation}
0\rightarrow M'_0\rightarrow\cdots \rightarrow M'_D\rightarrow 0
\label{M'}
\end{equation}
of free right $\cala$-modules with $M'_n\simeq E^\ast_n \otimes \cala$ where for any vector space $F$, one denotes by $F^\ast$ its dual vector space.\\

The algebra $\cala\in \mathbf{GrAlg}$ is said to be {\sl Gorenstein} whenever one has
\[
\left\{
\begin{array}{l}
H^n(M')=0, \> \> \> \text{for}\> \> \> n\not= D\\
H^D(M')=\mathbb K
\end{array}
\right.
\]
which reads $\Ext^n_\cala(\mathbb K, \cala)=\delta^{nD}\mathbb K$ by definition $(\delta^{nD}=0$ for $n\not=D$ and $\delta^{DD}=1$). This means that one has
\[
E^\ast_{D-n}\simeq E_n
\]
which is our version of the Poincaré duality.\\

Notice that one has
\[
\Ext^n_\cala(\mathbb K, \mathbb K)\simeq E^\ast_n
\]
which follows easily from the definitions.

\section{Quadratic algebras}

\subsection{Koszul duality and the Koszul complex}

A (homogeneous) {\sl quadratic algebra} $\cala$ is a graded algebra of the form
\[
\cala=A(E,R)=T(E)/[R]
\]
where $E$ is a finite-dimensional vector space and where $R\subset E\otimes E$ is a linear subspace of $E\otimes E$ ($[R]$ denotes as before the ideal of $T(E)$ generated by $R$), see for instances \cite{man:1988}, \cite{pol-pos:2005}.\\

Given a quadratic algebra $\cala=A(E,R)$, one defines another quadratic algebra $\cala^!$ called the {\sl Koszul dual algebra} of $\cala$ by setting
\[
\cala^!=A(E^\ast,R^\perp)
\]
with $R^\perp$ defined by
\[
R^\perp = \{\omega\in E^\ast \otimes E^\ast \vert \omega(r)=0,\>\> \> \forall r\in R\}
\]
where we have made the identification
$E^\ast\otimes E^\ast=(E\otimes E)^\ast$ which is allowed since one has $\dim(E)<\infty$.\\

One has $\cala^!=\otimes_{n\geq 0} \cala^!_n$ and the subspace $\cala^!_n$ of elements of degree $n$ of $\cala^!$ is given by
\[
\cala^!_n=E^{\ast\otimes^n}/\sum_k E^{\ast\otimes^k}\otimes R^\perp \otimes E^{\ast\otimes^{n-k-2}}
\]
 which is equivalent for its dual to 
\begin{equation}
\cala^{!\ast}_n=\cap_{0\leq k\leq n-2} E^{\otimes^k}\otimes R\otimes E^{\otimes^{n-k-2}}
\label{ddKn}
\end{equation}
for any $n\in \mathbb N$. The {\sl Koszul complex} of $\cala$ is then defined to be the chain complex $K(\cala)$ of free left $\cala$-module given by
\begin{equation}
\cdots \stackrel{d}{\rightarrow} \cala\otimes \cala^{!\ast}_{n+1}\stackrel{d}{\rightarrow} \cala\otimes \cala^{!\ast}_n \stackrel{d}{\rightarrow}\cdots
\label{Kc}
\end{equation}
where $d$ is induced by the left $\cala$-module homomorphisms
\[
d:\cala\otimes E^{\otimes^{n+1}}\rightarrow \cala\otimes E^{\otimes^n}
\]
given by setting
\[
d(a\otimes (e_0\otimes e_1\otimes \cdots \otimes e_n))=ae_0 \otimes (e_1\otimes \cdots \otimes e_n)
\]
for $n\in \mathbb N$ and $d(a)=0$ for $a\in \cala$.

\subsection{Koszul quadratic algebras}

A quadratic algebra $\cala$ is said to be {\sl Koszul} whenever its Koszul complex $K(\cala)$ is acyclic in positive degrees i.e. iff one has $H_n(K(\cala))=0$ for $n\geq 1$. If $\cala$ is Koszul, one has then $H_0(K(\cala))=\mathbb K$ so that the Koszul complex gives a free resolution
\[
K(A)\rightarrow \mathbb K\rightarrow 0
\]
of the trivial left module $\mathbb K$ which is in fact a minimal projective resolution of $\mathbb K$, \cite{pri:1970}. Thus, if $\cala$ is Koszul of finite global dimension $D$, one has $K_n(\cala)=0$ for $n>D$ and $K_D(\cala)\not=0$.\\

Let $\cala$ be a quadratic algebra and let us apply the functor $\Hom_\cala(\bullet, \cala)$ to the Koszul complex which is a chain complex of left $\cala$-modules. One then obtains the cochain complex $L(\cala)$ of right $\cala$-modules given by
\begin{equation}
\cdots \stackrel{\delta}{\rightarrow} \cala^!_n\otimes \cala\stackrel{\delta}{\rightarrow} \cala^!_{n+1}\otimes \cala\stackrel{\delta}{\rightarrow}\cdots
\label{Lc}
\end{equation}
where $\delta$ is the left multiplication by $x_\lambda^\ast \otimes x^\lambda$ for $(x^\lambda)$  a basis of $E=\cala_1$ with dual basis $(x^\ast_\lambda)$ of $E^\ast=\cala^!_1$.\\

Assume now that $\cala$ is a Koszul algebra of finite global dimension $D$. Then it follows from above that $\cala$ is Gorenstein if and only if one has
\[
H^n(L(\cala))=\delta^{nD}\mathbb K
\]
for the cohomology of $L(\cala)$. Thus for Koszul algebras of finite global dimension, the Poincaré duality property is controlled by the cohomology of $L(\cala)$.

\subsection{Twisted potentials for Koszul-Gorenstein algebras}

Let $V$ be a vector space and let $n\geq 1$ be a positive integer, then a\linebreak[4] $(n+1)$-linear form $w$ on $V$ is said to be {\sl preregular} \cite{mdv:2005}, \cite{mdv:2007} iff it satisfies the following conditions (i) and (ii).\\
(i) If $X\in V$ is such that $w(X,X_1,\dots, X_n)=0$ for any $X_1,\dots,X_n\in V$, then $X=0$.\\
(ii) There is an element $Q_w\in GL(V)$ such that one has
\[
w(X_0,\dots,X_{n-1},X_n)=w(Q_wX_n,X_0,\dots,X_{n-1})
\]
for any $X_0,\dots, X_n\in V$.\\

It follows from (i) that $Q_w$ as in (ii) is unique. Property (i) will be refered to as 1-site nondegeneracy while (ii) will be refered to as twisted cyclicity.\\

Let $E$ be a finite-dimensional vector space, then an element $w$ of $E^{\otimes^D}$ is the same thing as a $D$-linear form on the dual $E^\ast$ of $E$. To make contact with the terminology of \cite{gin:2006}we will say that $w$ is a {\sl twisted potential} of degree $D$ on $E$ if the corresponding $D$-linear form on $E^\ast$ is preregular.\\

Let $w\in E^{\otimes^D}$ be a twisted potential and let $w_{\lambda_1\dots \lambda_D}$ be its components in the basis $(x^\lambda)_{\lambda\in \{1,\dots,d\}}$ of $E$, i.e. one has $w=w_{\lambda_1\dots \lambda_D}x^{\lambda_1}\otimes \dots \otimes x^{\lambda_D}$. Let $R_w\subset E\otimes E$ be the subspace of $E\otimes E$ spanned by the  $w_{\lambda_1\dots \lambda_{D-2}\mu\nu} x^\mu\otimes x^\nu$ for $\lambda_1,\dots,\lambda_{D-2}\in \{1,\dots,d\}$. In this way one can associate to $w$ the quadratic algebra $\cala(w,2)=A(E,R_w)$. The content of the main result of \cite{mdv:2005} or of \cite{mdv:2007} (Theorem 4.3 of \cite{mdv:2005} or Theorem 11 of \cite{mdv:2007}) applied to quadratic algebras (case $N=2$) is the following

\begin{theorem}\label{POT}
Let $\cala$ be a quadratic Koszul algebra of finite global dimension $D$ which is Gorenstein. Then $\cala=\cala(w,2)$ for some twisted potential $w$ of degree $D$.
\end{theorem}

In the case where $Q_w=(-1)^{D+1}$, $\cala$ is a quadratic Calabi-Yau algebra.\\

This result is as mentioned the particular case for $N=2$ of a general result for $N$-homogeneous algebras (i.e. relations of degree N), see in \cite{mdv:2005}, \cite{mdv:2007}. The case $N=3$ (cubic algebras) contains the important example of the Yang-Mills algebra \cite{ac-mdv:2002b} which is a graded Calabi-Yau algebra in the sense of \cite{gin:2006}.\\

As pointed out in \cite{mdv:2007}, $\bbbone \otimes w\>\> (\bbbone\in \cala)$ has the interpretation of a (twisted) noncommutative volume form, see Theorem 10 of 
\cite{mdv:2007}.\\

Thus the Poincaré duality corresponding to the Gorenstein property implies for the Koszul algebras that they are derived from twisted potentials. Let us now give some examples

\subsection{Examples}

\begin{enumerate}
\item {\sl The symmetric algebra} $\cala=SE$. One has $\cala=SE=A(E,R)$ with $R=\wedge^2 E\subset E\otimes E$, therefore $\cala^!=\wedge E^\ast$ is the exterior algebra of $E^\ast$. It follows that the Koszul complex is $SE\otimes \wedge^\bullet E$ with the Koszul differential $d:SE\otimes \wedge^{n+1} E\rightarrow SE\otimes \wedge^n E$. But $SE\otimes \wedge^\bullet E$ is also the algebra of polynomial differential forms and, if $\delta$ denotes the exterior differential, the derivation $d\delta + \delta d$ of degree 0 coincides with the total degree i.e. $(d\delta+\delta d)x=(r+s)x$ for $x\in S^r E\otimes \wedge^s E$. Thus $\delta$ is a homotopy for $d$ and $d$ is a homotopy for $\delta$ whenever $r+s\not=0$. This implies that $SE$ is Koszul of global dimension $D=\dim(E)$ and this also implies the formal Poincaré lemma.\\

The algebra $\cala=SE$ has also the Gorenstein property which here reduces to the usual Poincaré duality property. The corresponding (twisted) potential reads 
\[
w=\varepsilon_{\lambda_1\dots \lambda_D}x^{\lambda_1}\otimes \dots \otimes x^{\lambda_D}
\]
where $D=\dim(E)$ and where $\varepsilon_{\lambda_1\dots \lambda_D}$ is completely antisymmetric with $\varepsilon_{1\dots D}=1$. It is clear that $v=\bbbone\otimes w$ is a volume form in the classical sense.

\item {\sl The tensor algebra} $\cala=T(E)$. One has $\cala=T(E)=A(E,0)$ so that $\cala^!=\mathbb K\bbbone \oplus E=A(E,E\otimes E)$ with trivial product between the elements of $E$. The Koszul complex is
\[
0\rightarrow T(E)\otimes E \stackrel{m}{\rightarrow} T(E)\rightarrow 0
\]
and it is obvious that $m$ is injective so $T(E)$ is Koszul of global dimension $D=1$ but it is clearly not Gorenstein, i.e. one has not the Poincaré duality here.

\item {\sl Koszul duals of Koszul algebras}. It is not hard to show that if $\cala$ is a quadratic algebra, then $\cala$ is Koszul if and only if $\cala^!$ is Koszul. Thus for instance $\wedge E$ and $\mathbb K\bbbone \oplus E$ are Koszul; however they are not Gorenstein, (they are not of finite global dimension). It is worth noticing here that this property is specific for quadratic algebras. Indeed there is a notion of Koszulity for $N$-homogeneous algebras which was introduced in \cite{ber:2001a} and a generalization of the Koszul duality for these algebras defined in \cite{ber-mdv-wam:2003} but the Koszulity is not stable by Koszul duality for $N\geq 3$.

\item {\sl A deformed symmetric algebra}. Let $\cala$ be the algebra generated by the 3 elements $\nabla_0,\nabla_1,\nabla_2$ with the relations
\begin{equation}
\left\{
\begin{array}{l}
\mu^2\nabla_2\nabla_0-\nabla_0\nabla_2=0\\
\mu^4\nabla_1\nabla_0-\nabla_0\nabla_1=0
\\
\mu^4\nabla_2\nabla_1-\nabla_1\nabla_2=0
\end{array}
\right.
\label{hrW}
\end{equation}
Where $\mu\in \mathbb K$ with $\mu\not=0$. This is a quadratic algebra with Koszul dual $\cala^!$ which is generated by the 3 elements $\omega_0$, $\omega_1,\omega_2$ with the relations
\begin{equation}
\left\{
\begin{array}{l}
\omega^2_0=0, \omega^2_1=0, \omega^2_2=0\\
\omega_2\omega_0 + \mu^2\omega_0\omega_2=0\\
\omega_1\omega_0+\mu^4\omega_0\omega_1=0\\
\omega_2\omega_1+\mu^4\omega_1\omega_2=0
\end{array}
\right.
\label{dhrW}
\end{equation}
It can be shown that $\cala$ is Koszul of global dimension 3 and is Gorenstein\cite{gur:1990},\cite{wam:1993}, \cite{mdv:2010}. The corresponding (twisted) potential reads
\[
\begin{array}{lll}
w=\mu^2(\nabla_1\otimes \nabla_2\otimes \nabla_0 & + &\nabla_2\otimes \nabla_0\otimes \nabla_1+\mu^{-6}\nabla_0\otimes \nabla_1\otimes \nabla_2)\\
& - &(\nabla_0\otimes \nabla_2\otimes \nabla_1+\nabla_1\otimes \nabla_0\otimes \nabla_2+\mu^6\nabla_2\otimes \nabla_1\otimes \nabla_0)
\end{array}
\]
while $Q_w\in GL(3, \mathbb K)$ is given by the diagonal matrix with
\[
(Q_w)^0_0=\mu^{-6},\>\> (Q_w)^1_1=1,\>\> (Q_w)^2_2=\mu^6
\]
in the basis $(\nabla_0,\nabla_1,\nabla_2)$.
\end{enumerate}

\section{Nonhomogeneous quadratic algebras}

\subsection{Poincaré-Birkhoff-Witt (PBW) property}

A {\sl nonhomogeneous quadratic algebra} \cite{pos:1993}, \cite{bra-gai:1996}, 
\cite{flo:2006}is an algebra $\fraca$ of the form
\[
\fraca=A(E,P)=T(E)/[P]
\]
where $E$ is a finite dimensional vector space and where $P\subset F^2(T(E))$ is a linear subspace of $F^2(T(E))=\oplus^2_{m=0} E^{\otimes^m}$. Here and in the following the tensor algebra $T(E)$ is endowed with its natural filtration $F^n(T(E))=\oplus_{m\leq n} E^{\otimes^m}$ associated to its graduation. The filtration of $T(E)$ induces a filtration $F^n(\fraca)$ of $\fraca$ and the graded algebra
\[
\gr (\fraca)=\oplus_n F^n(\fraca)/F^{n-1}(\fraca)
\]
is the {\sl associated graded algebra} to the filtered algebra $\fraca$. Let $R\subset E \otimes E$ be the image of $P$ by the canonical projection of $F^2(T(E))$ onto $E\otimes E$. Then $\cala=A(E,R)=T(E)/[R]$ is a (homogeneous) quadratic algebra which will be called the {\sl quadratic part} of $\fraca$. There is a canonical homomorphism
\[
can:\cala\rightarrow \gr(\fraca)
\]
of graded algebra which is surjective. The algebra is said to have the {\sl Poincaré-Birkhoff-Witt (PBW) property} if this canoncal homomorphism is injective, i.e. whenever $can$ is an isomorphism. One has the following theorem 
\cite{bra-gai:1996}.

\begin{theorem}\label{PBW}
Let $\fraca$ and $\cala$ be as above. If $\fraca$ has the PBW property then the following condition $\mathrm{(i)}$ and $\mathrm{(ii)}$ are satisfied :\\
$\mathrm{(i)}$ \hspace{0,1cm} $P\cap F^1(T(E))=0$\\
$\mathrm{(ii)}$ \hspace{0,1cm} $(PE+EP)\cap F^2(T(E))\subset P$.\\
If $\cala$ is Koszul then, conversely, conditions $\mathrm{(i)}$ and $\mathrm{(ii)}$ imply that $\fraca$ has the PBW property.
\end{theorem}
Condition (i) means that one has linear mappings $\varphi:R\rightarrow E$ and\linebreak[4] $\varphi_0 : R\rightarrow \mathbb K$ such that
\[
P=\{r-\varphi(r)-\varphi_0(r)\bbbone\>\vert\> r\in R\}
\]
i.e. $P$ is obtained by adding second members of lower degrees to the quadratic relations $R$. Concerning condition (ii) one has the following.

\begin{proposition}\label{II}
Assume that Condition $\mathrm{(i)}$ of the last theorem is satisfied. Then, Condition $\mathrm{(ii)}$ above is equivalent to the following conditions
$\mathrm{(a)}$, $\mathrm{(b)}$ and $\mathrm{(c)}$:\\
$\mathrm{(a)}$\hspace{0,1cm} $(\varphi\otimes I-I\otimes \varphi)(R\otimes E\cap E\otimes R)\subset R$\\
$\mathrm{(b)}$ \hspace{0,1cm} $(\varphi\circ (\varphi\otimes I-I\otimes \varphi)+ (\varphi_0\otimes I-I\otimes \varphi_0))(R\otimes E\cap E\otimes R)=0$\\
$\mathrm{(c)}$ \hspace{0,1cm} $\varphi_0\circ (\varphi\otimes I-I\otimes \varphi) (R\otimes E\cap E\otimes R)=0$
where $I$ denotes the identity mapping of $E$ onto itself.
\end{proposition}
A nonhomogeneous quadratic algebra $\fraca$ with quadratic part $\cala$ is said to be {\sl Koszul} if $\fraca$ has the PBW property and $\cala$ is Koszul.

\subsection{Nonhomogeneous Koszul duality}

Let $\fraca=A(E,P)$ be a nonhomogeneous quadratic algebra with quadratic part $\cala=A(E,R)$, assume that Condition (i) of Theorem \ref{PBW} is satisfied and let $\varphi:R\rightarrow E$ and $\varphi_0:R\rightarrow \mathbb K$ be as in 3.1. Consider the transposed $\varphi^t:E^\ast\rightarrow R^\ast$ and $\varphi^t_0:\mathbb K\rightarrow R^\ast$ of $\varphi$ and $\varphi_0$ and notice that one has by definition  of $\cala^!$ that $\cala^!_1=E^\ast$, $\cala^!_2=R^\ast$ and $\cala^!_3=(R\otimes E\cap E\otimes R)^\ast$, (see in 2.1 formula (\ref{ddKn})), so one can write (the minus sign is put here to match the usual conventions)
\begin{equation}
-\varphi^t:\cala^!_1\rightarrow \cala^!_2,\> -\varphi^t_0(1)=F\in \cala^!_2
\label{Dual}
\end{equation}
and one has the following result \cite{pos:1993}.
\begin{theorem}\label{KD}
Conditions $\mathrm{(a)}$, $\mathrm{(b)}$ and $\mathrm{(c)}$ of Proposition \ref{II} are equivalent to the following conditions $\mathrm{(a')}$, $\mathrm{(b')}$ and $\mathrm{(c')}$ :\\
$\mathrm{(a')}$\hspace{0,1cm} $-\varphi^t$ extends as an antiderivation $\delta$ of $\cala^!$\\
$\mathrm{(b')}$\hspace{0,1cm} $\delta^2(x)=[F,x],\>\> \forall x\in \cala^!$\\
$\mathrm{(c')}$\hspace{0,1cm} $\delta(F)=0$.
\end{theorem}

A graded algebra equipped with an antiderivation $\delta$ of degree 1 and an element $F$ of degree 2 satisfying the conditions $\mathrm{(b')}$ and $\mathrm{(c')}$ above is refered to as a {\sl curved graded differential algebra} \cite{pos:1993}.\\

Thus Theorem \ref{KD} combined with Theorem \ref{PBW} and Proposition \ref{II} means that the correspondance $\fraca\mapsto (\cala^!,\delta,F)$ define a contravariant functor from the category of nonhomogeneous quadratic algebras satisfying the conditions (i) and (ii) of Theorem \ref{PBW} to the category of curved differential quadratic algebras (with the obvious appropriate notions of morphism). One can summarize the Koszul duality of \cite{pos:1993} for non homogeneous quadratic algebras by the following.

\begin{theorem}\label{POS}
The above correspondence defines an anti-isomorphism between the category of nonhomogeneous quadratic algebras satisfying Conditions $\mathrm{(i)}$ and $\mathrm{(ii)}$ of Theorem \ref{PBW} and the category of curved differential quadratic algebras which induces an anti-isomorphism between the category of nonhomogeneous quadratic Koszul algebras and the category of curved differential quadratic Koszul algebras.
\end{theorem}

There are  two important classes of nonhomogeneous quadratic algebras $\fraca$ satisfying the conditions (i) and (ii) of Theorem \ref{PBW}. The first one corresponds to the case $\varphi_0=0$ which is equivalent to $F=0$ while the second one corresponds to $\varphi=0$ which is equivalent to $\delta=0$. An algebra $\fraca$ of the first class is called a {\sl quadratic-linear algebra} 
\cite{pol-pos:2005} and corresponds to a differential quadratic algebra $(\cala^!,\delta)$ while an algebra $\fraca$ of the second class corresponds to
a quadratic algebra $\cala^!$ equipped with a central element $F$ of degree 2.

\subsection{Examples}

\begin{enumerate}
\item {\sl Universal enveloping algebras of Lie algebras}.
Let $\fracg$ be a finite-dimen\-sional Lie algebras then its universal enveloping algebra $\fraca=U(\fracg)$ is Koszul quadratic-linear. Indeed one has $\cala=S\fracg$ which is a Koszul quadratic algebra while the PBW property is here the classical PBW property of $U(\fracg)$. The corresponding differential quadratic algebra $(\cala^!,\delta)$ is $(\wedge\fracg^\ast,\delta)$, i.e. the exterior algebra of the dual vector space $\fracg^\ast$ of $\fracg$ endowed with the Koszul differential $\delta$. Notice that this latter differential algebra is the basic building block to construct the Chevalley-Eilenberg cochain complexes. Notice also that $\cala=S\fracg$ is not only Koszul but is also Gorenstein (Poincaré duality property).

\item {\sl Adjoining a unit element to an associative algebra}. Let $A$ be a finite-dimensional associative algebra and let 
\[
\fraca=\tilde A=T(A) / \left[\{ x \otimes y-xy, y \in A \} \right]
\]
 be the algebra obtained by adjoining a unit $\bbbone$ to $A$ ($\tilde A=\mathbb K \bbbone\oplus A$, etc.). 
This is again a Koszul quadratic-linear algebra. Indeed the PBW property is here equivalent to the associativity of $A$ while the quadratic part is $\cala=T(A^\ast)^!$ which is again $\mathbb K\bbbone\oplus A$ as vector space but with a vanishing product between the elements of $A$ and is a Koszul quadratic algebra. The corresponding differential quadratic algebra $(\cala^!,\delta)$ is $(T(A^\ast),\delta)$ where $\delta$ is the antiderivation extension of minus the transposed $m^t:A^\ast\rightarrow A^\ast\otimes A^\ast$ of the product $m$ of $A$. Again $(T_+(A^\ast),\delta)$ is the basic building block to construct the Hochschild cochain complexes. Notice however that $\cala=T(A^\ast)^!$ is not Gorenstein (no Poincaré duality).

\item {\sl A deformed universal enveloping algebra}. Let $\fraca$ be the algebra generated by the 3 elements $\nabla_0,\nabla_1,\nabla_2$ with relations
\begin{equation}
\left \{
\begin{array}{l}
\mu^2\nabla_2\nabla_0-\nabla_0\nabla_2=\mu\nabla_1\\
\mu^4\nabla_1\nabla_0-\nabla_0\nabla_1=\mu^2(1+\mu^2)\nabla_0\\
\mu^4\nabla_2\nabla_1-\nabla_1\nabla_2=\mu^2(1+\mu^2)\nabla_2
\end{array}
\right.
\label{irW}
\end{equation}
This is again a Koszul quadratic-linear algebra with homogeneous part $\cala$ given by the example 4 of 2.4 which is Koszul-Gorenstein. The corresponding differential quadratic algebra $(\cala^!,\delta)$ is the algebra $\cala^!$ of Example 4 of 2.4 with quadratic relations (\ref{dhrW}) endowed with the differential $\delta$ given by
\begin{equation}
\left \{
\begin{array}{l}
\delta\omega_0+\mu^2(1+\mu^2)\omega_0\omega_1=0\\
\delta\omega_1+\mu \omega_0\omega_2=0\\
\delta\omega_2+\mu^2(1+\mu^2)\omega_1\omega_2=0
\end{array}
\right.
\label{diffW}
\end{equation}
which corresponds to the left covariant differential calculus on the twisted $SU(2)$ group of \cite{wor:1987b}.

\item {\sl Canonical commutation relations algebra}. Let $E=\mathbb K^{2n}$ with basis $(q^\lambda, p_\mu)$, $\lambda, \mu\in \{1,\dots,n\}$ and let $i\hbar\in \mathbb K$ with $i\hbar\not=0$. Consider the nonhomogeneous quadratic algebra
$\fraca$ generated by the $q^\lambda,p_\mu$ with relations
\[
q^\lambda q^\mu-q^\mu q^\lambda=0,\>\>p_\lambda p_\mu-p_\mu p_\lambda=0,\>\>q^\lambda p_\mu-p_\mu q^\lambda=i\hbar \delta^\lambda_\mu\bbbone
\]
for $\lambda,\mu\in \{1,\dots, n\}$. The quadratic part of $\fraca$ is the symmetric algebra $\cala=SE$ which is Koszul, the property (i) of Theorem \ref{PBW} is obvious and one has $\varphi=0$ and $\varphi_0$ is such that its transposed $\varphi^t_0$ is given by
\[
-\varphi^t_0(1)=F=-(i\hbar)^{-1} q^\ast_\lambda\wedge  p^{\lambda^\ast}
\]
which is central in $\cala^!=\wedge(E^\ast)$ where $(q^\ast_\lambda,p^{\mu\ast})$ is the dual basis of $(q^\lambda,p_\mu)$. This implies that $\fraca$ has the PBW property and therefore is Koszul.

\item {\sl Clifford algebra (C.A.R. algebra)}. Let $E=\mathbb K^n$ with canonical basis $(\gamma_\lambda)$, $\lambda\in \{1,\dots,n\}$ and consider the nonhomogeneous quadratic algebra $\fraca=C(n)$ generated by the elements $\gamma_\lambda$, $\lambda\in \{1,\dots, n\}$ with relations 
\[
\gamma_\mu \gamma_\nu+\gamma_\nu \gamma_\mu=2\delta_{\mu\nu}\bbbone
\]
for $\mu,\nu\in\{1,\dots,n\}$. The quadratic part of $\fraca$ is then the exterior algebra $\cala=\wedge E$ which is Koszul, the property (i) of Theorem \ref{PBW} is obvious and one has again $\varphi=0$ and $\varphi^t_0$ is given by
\[
-\varphi^t_0(1)=F=-\frac{1}{2}\sum \gamma^{\lambda\ast} \vee \gamma^{\lambda\ast}
\]
which is a central element of $\cala^!=SE^\ast$ (which is commutative). It again follows that $\fraca$ is Koszul (i.e. PBW + $\cala$ Koszul).

\item {\sl Remarks on the generic case}. Let $\cala$ be a (homogeneous) quadratic algebra which is Koszul. In general (for generic $\cala$) any nonhomogeneous quadratic algebra $\fraca$ which has $\cala$ as quadratic part and has the PBW property is such that one has both $\varphi\not=0$ and $\varphi_0\not=0$ or is trivial in the sense that it coincides with $\cala$, i.e. $\varphi=0$ and $\varphi_0=0$. This is the case for instance when $\cala$ is the 4-dimensional Sklyanin algebra \cite{skl:1982}, \cite{smi-sta:1992} for generic values of its parameters \cite{bel-ac-mdv:2011}.

Thus, Examples 1, 2, 3, 4, 5 above are rather particular from this point of view. However the next section will be devoted to a generalization of Lie algebra which has been introduced in \cite{mdv-lan:2011} and which involves quadratic-linear algebras, i.e. for which $\varphi_0=0$.

\end{enumerate}

\section{Lie prealgebras}

\subsection{Prealgebras} By a (finite-dimensional) {\sl prealgebra} we here mean a triple $(E,R,\varphi)$ where $E$ is a finite-dimensional vector space, $R\subset E\otimes E$ is a linear subspace of $E\otimes E$ and $\varphi:R\rightarrow E$ is a linear mapping of $R$ into $E$. Given a supplementary $R'$ to $R$ in $E\otimes E$, $R\oplus R'=E\otimes E$, the corresponding projector $P$ of $E\otimes E$ onto $R$ allows to define a bilinear product $\varphi\circ P:E\otimes E\rightarrow E$, i.e. a structure of algebra on $E$. The point is that there is generally no natural supplementary of $R$. Exception are $R=E\otimes E$ of course and $R=\wedge^2 E\subset E\otimes E$ for which there is the canonical $GL(E)$-invariant supplementary $R'=S^2E\subset E\otimes E$ which leads to an antisymmetric product on $E$, (e.g. case of the Lie algebras). \\

Given a prealgebra $(E,R,\varphi)$, there are two natural associated algebras :
\begin{enumerate}
\item
The nonhomogeneous quadratic algebra
\[
\fraca_E=T(E)/[\{r-\varphi(r)\>\> \vert\>\> r\in R\}]
\]
which will be called its {\sl enveloping algebra}.

\item
The quadratic part $\cala_E$ of $\fraca_E$
\[
\cala_E=T(E)/[R],
\]
where the prealgebra $(E,R,\varphi$) is also simply denoted by $E$ when no confusion arises.

\end{enumerate}

The enveloping algebra $\fraca_E$ is a filtered algebras as explained before but it is also an augmented algebra with augmentation
\[
\varepsilon:\fraca_E\rightarrow \mathbb K
\]
induced by the canonical projection of $T(E)$ onto $T^0(E)=\mathbb K$. One has the surjective homomorphism
\[
can:\cala_E\rightarrow \gr(\fraca_E)
\]
of graded algebras.\\

In the following we shall be mainly interested on prealgebras such that their enveloping algebras are quadratic-linear, i.e. satisfy Condition (ii) of Theorem \ref{PBW}, (Condition (i) being satisfied by construction). If $(E,R,\varphi)$ is such a prealgebra, to $\fraca_E$ corresponds the differential quadratic algebra $(\cala^!_E,\delta)$ (as in Section 3) where $\delta$ is the antiderivation extension of minus the transposed $\varphi^t$ of $\varphi$.\\

Notice that if $\fraca_E$ has the PBW property one has
\[
E=F^1(\fraca_E)\cap \ker(\varepsilon)
\]
so that the canonical mapping of the prealgebra $E$ into its enveloping algebra $\fraca_E$ is then an injection.

\subsection{Lie prealgebras}
A prealgebra $(E,R,\varphi)$ will be called a {\sl Lie prealgebra} \cite{mdv-lan:2011} if the following conditions (1) and (2) are satisfied :\\
(1) The quadratic algebra $\cala_E=A(E,R)$ is Koszul of finite global dimension and is Gorenstein (Poincaré duality). \\
(2) The enveloping algebra $\fraca_E$ has the PBW property.\\

If $E=(E,R,\varphi)$ is a Lie prealgebra then $\fraca_E$ is a Koszul quadratic linear algebra, so to $(E,R,\varphi)$ one can associate the differential quadratic algebra $(\cala^!_E,\delta)$ and one has the following theorem \cite{mdv-lan:2011}:
\begin{theorem}\label{LPD}
The correspondence $(E,R,\varphi)\mapsto (\cala^!_E,\delta)$ defines an anti-isomorphism between the category of Lie prealgebra and the category of differential quadratic Koszul Frobenius algebras.
\end{theorem}
This is a direct consequence of Theorem \ref{POS} and of the Koszul Gorenstein property of $\cala_E$ by using \cite{smi:1996}.\\

Let us remind that a {\sl Frobenius algebra} is a finite-dimensional algebra $\cala$ such that as left $\cala$-modules $\cala$ and its vector space dual $\cala^\ast$ are isomorphic (the left $\cala$-module structure of $\cala^\ast$ being induced by the right $\cala$-module structure of $\cala$). Concerning the graded connected case one has the following classical useful result.

\begin{proposition}\label{GF}
Let $\cala=\oplus_{m\geq 0} \cala_m$ be a finite-dimensional graded connected algebra with $\cala_D\not=0$ and $\cala_n=0$ for $n>D$. Then the following conditions $\mathrm{(i)}$ and $\mathrm{(ii)}$ are equivalent :\\
$\mathrm{(i)}$ $\cala$ is Frobenius,\\
$\mathrm{(ii)}$ $\dim(\cala_D)=1$  and $(x,y)\mapsto (xy)_D$ is nondegenerate, where $(z)_D$ denotes the component on $\cala_D$ of $z\in \cala$.
\end{proposition}

\subsection{Examples and counterexamples}

\begin{enumerate}

\item {\sl Lie algebras.} It is clear that a Lie algebra $\fracg$ is canonically a Lie prealgebra $(\fracg, R,\varphi)$ with $R=\wedge^2\fracg\subset \fracg\otimes \fracg,\> \varphi=[\bullet, \bullet]$, $\fraca_\fracg=U(\fracg)$ and $\cala_\fracg=S\fracg$,(see Example 1 in 3.3).

\item {\sl Associative algebras are not Lie prealgebras.} An associative algebra $A$ is clearly a prealgebra $(A,A\otimes A,m)$ with enveloping algebra $\fraca_A=\tilde A$ as in Example 2 of 3.3 but $\cala_A=T(A^\ast)^!=\mathbb K\bbbone\oplus A$ is not Gorenstein although it is Koszul as well as $\fraca_A=\tilde A$, (see the discussion of Example 2 in 3.3). The missing item is here the Poincaré duality.

\item {\sl A deformed version of Lie algebras.} The algebra $\fraca$ of Example 3 of 3.3 is the enveloping algebra of a Lie prealgebra $(E,R,\varphi)$ with $E=\mathbb K^3$, $R\subset E\otimes E$ generated by\\
$r_1=\mu^2\nabla_2\otimes \nabla_0-\nabla_0\otimes \nabla_2$\\
$r_0=\mu^4\nabla_1\otimes \nabla_0-\nabla_0\otimes \nabla_1$\\
$r_2=\mu^4\nabla_2\otimes \nabla_1-\nabla_1\otimes \nabla_2$\\
and $\varphi$ given by
\[
\varphi(r_1)=\mu\nabla_1,\>\>\>\> \varphi(r_0)=\mu^2(1+\mu^2)\nabla_0,\>\> \>\>\varphi(r_2)=\mu^2(1+\mu^2)\nabla_2
\]
where $(\nabla_0,\nabla_1,\nabla_2)$ is the canonical basis of $E$.

\item {\sl Differential calculi on quantum groups.} More generally most differential calculi on the quantum groups can be obtained via the duality of Theorem 6 from Lie prealgebras. In fact the Frobenius property is generally straightforward to verify, what is less obvious to prove is the Koszul property.
\end{enumerate}

\subsection{Generalization of Chevalley-Eilenberg complexes}

Throughout this subsection, $E=(E,R,\varphi)$ is a fixed Lie prealgebra, its enveloping algebra is simply denoted by $\fraca$ with quadratic part denoted by $\cala$ and the associated differential quadratic Koszul Frobenius algebra is $(\cala^!,\delta)$.\\

A {\sl left representation} of the Lie prealgebra $E=(E,R,\varphi)$ is a left $\fraca$-module. Let $V$ be a left representation of $E=(E,R,\varphi)$, let $(x^\lambda)$ be a basis of $E$ with dual basis $(\omega_\lambda)$ of $E^\ast=\cala^!_1$. One has
\[
x^\mu x^\nu \Phi \otimes \omega_\mu\omega_\nu+ x^\lambda \Phi\otimes \delta \omega_\lambda=0
\]
for any $\Phi\in V$. This implies that one defines a differential of degree 1 on $V\otimes \cala^!$ by setting
\[
\delta_V(\Phi\otimes \alpha)=x^\lambda \Phi \otimes \omega_\lambda\alpha + \Phi \otimes \delta \alpha
\]
so $(V\otimes \cala^!,\delta_V)$ is a cochain complex. These cochain complexes generalize the Chevalley-Eilenberg cochain complexes. Given a right representation of $E$, that is a right $\fraca$-module $W$, one defines similarily the chain complex $(W\otimes \cala^{!\ast},\delta_W)$, remembering that $\cala^{!\ast}$ is a graded coalgebra.\\

One has the isomorphisms
\[
\left\{
\begin{array}{l}
H^\bullet (V\otimes \cala^!) \simeq \Ext^\bullet_\fraca(\mathbb K, V)\\
H_\bullet(W\otimes \cala^{!\ast})\simeq \Tor^\fraca_\bullet (W,\mathbb K)
\end{array}
\right.
\]
which implies that one has the same relation with the Hochschild cohomology and the Hochschild homology of $\fraca$ as the relation of the (co-)homology of a Lie algebra with the Hochschild (co-)homology of its universal enveloping algebra.

\section{Conclusion}

As pointed out before many results of this article, in particular all the results of Section 2 generalize to the case of relations of degree $N\geq 2$ (instead of 2) 
\cite{mdv:2007}. Furthermore the results of \cite{mdv:2007} have been generalized in \cite{boc-sch-wem:2010} to the quiver case. In fact the analysis of \cite{ber-gin:2006} suggests that one can generalize most points by replacing the ground field $\mathbb K$ by a von Neumann regular ring $\mathbf R$ and replacing the tensor algebras of vector spaces by tensor algebra over $\mathbf R$ of $\mathbf R$-bimodules. The case of quiver being $\mathbf R=\mathbb K^{vertex}$.\\

A ring $\mathbf R$ is said to be {\sl von Neumann regular} whenever for any $a\in \mathbf R$ there is an $x\in \mathbf R$ such that $axa=a$. Semisimple rings are von Neumann regular. An infinite product of fields and the endomorphism ring $\End(E)$ of an infinite-dimensional vector space $E$ are examples of von Neumann regular rings which are not semisimple rings \cite{wei:1994}, \cite{ber-gin:2006}.\\

Concerning the nonhomogeneous case, it is worth noticing here that for nonhomogeneous relations of degree $N\geq 3$ one has no satisfying generalization of the Koszul duality of Positselski for the moment.

\newpage

\end{document}